\documentclass[a4paper,11pt,oneside]{amsart}
\usepackage{amsmath,amsfonts,amssymb,amsthm,amscd,latexsym,euscript,xypic}
\usepackage[all]{xy}
\newtheorem{theorem}{Theorem}[section]
\newtheorem{corollary}[theorem]{Corollary}
\newtheorem{definition}[theorem]{Definition}
\newtheorem{example}[theorem]{Example}
\newtheorem{proposition}[theorem]{Proposition}
\newtheorem{remark}[theorem]{Remark}
\newcommand{\Ae}{{\mathcal A}}

\newcommand{\Ce}{{\mathcal C}}
\newcommand{\De}{{\mathcal D}}
\newcommand{\Ee}{{\mathcal E}}
\newcommand{\Fe}{{\mathcal F}}
\newcommand{\He}{{\mathcal H}}
\newcommand{\Le}{{\mathcal L}}

\newcommand{\Xe}{{\mathcal X}}
\newcommand{\Hom}{{\rm Hom}}
\newcommand{\OR}{{\rm Ord}}
\newcommand{\rank}{{\rm rank}}
\newcommand{\identity}{{\rm id}}
\newcommand{\Gra}{{\bf Gra}}
\newcommand{\TC}{{\rm TC}}
\newcommand{\Sets}{{\bf Sets}}

\newcommand{\Z}{{\mathbb Z}}
\newcommand{\Id}{{\rm Id}}

\numberwithin{equation}{section}

\begin{document}

\title{Epireflections and supercompact cardinals}
\author{Joan Bagaria, Carles Casacuberta, and Adrian R. D. Mathias}
\thanks{The authors were supported by the Spanish 
Ministry of Education and Science under MEC-FEDER grants MTM2004-03629 and MTM2005-01025,
and by the Generalitat de Catalunya under grants 2005SGR-00606 and 2005SGR-00738.
The hospitality of the Centre de Recerca Matem\`atica (CRM) is also acknowledged.}

\begin{abstract}
We prove that, under suitable assumptions on a category~$\Ce$,
the existence of supercompact cardinals implies that every 
absolute epireflective class of objects of $\Ce$ is a small-orthogonality class.
More precisely, if $L$ is a localization functor on an accessible category 
$\Ce$ such that the unit morphism $X\to LX$
is an extremal epimorphism for all~$X$, and the class of $L$-local objects is
defined by an absolute formula with parameters, then the existence of a supercompact
cardinal above the cardinalities of the parameters implies that $L$ is a localization 
with \hbox{respect} to some set of morphisms.
\end{abstract}

\maketitle
\section{Introduction}

The answers to certain questions in infinite abelian group theory
are known to depend on set theory. For example, the question whether torsion
theories are necessarily singly generated or singly cogenerated was
discussed in~\cite{DG}, where the existence or nonexistence of
measurable cardinals played a significant role. 
In a different direction, conditions under which cotorsion pairs
are generated or cogenerated by a set were studied in~\cite{EST}.
Other algebraic problems whose answer involves set-theoretical assumptions
can be found in~\cite{EM}.

In homotopy theory, it was asked around 1990 if every functor
on simplicial sets which is idempotent up to homotopy is equivalent to
$f$-localization for some map $f$ (see \cite{FarjounBCAT} and
\cite{Farjounbook} for terminology and details). Although this may not
seem a set-theoretical question, the following counterexample was given in~\cite{CSS}:
Under the assumption that measurable cardinals do not exist,
the functor $L$ defined as $LX=NP_{\Ae}(\pi X)$, where
$\pi$ denotes the fundamental groupoid, $N$ denotes the nerve, 
and $P_{\Ae}$ denotes reduction with respect to the proper
class $\Ae$ of groups of the form $\Z^{\kappa}/\Z^{<\kappa}$ for all
cardinals~$\kappa$, is not equivalent to localization with respect to
any set of maps.

The statement that measurable cardinals do not exist is consistent with the
Zermelo--Fraenkel axioms with the axiom of choice (ZFC), provided of course
that ZFC is itself consistent. However, many large-cardinal
assumptions, such as the existence of measurable cardinals, or
bigger cardinals, are used in mathematical practice, leading to
useful developments. Specifically, Vop\v{e}nka's principle
\cite{Jech2} implies that every homotopy idempotent functor on simplicial sets 
is an $f$-localization for some map~$f$; see \cite{CSS} for a proof of this claim.
Vop\v{e}nka's principle (one of whose forms is the statement that between the
members of every proper class of graphs there is at least one nonidentity map) 
has many other similar consequences, such as the fact that
all reflective classes in locally presentable categories are
small-orthogonality classes (i.e., orthogonal to some set of morphisms) \cite{AR}, 
or that all colocalizing subcategories
of triangulated categories derived from locally presentable Quillen model categories 
are reflective~\cite{CGR}.

In this article, we show that the existence of supercompact cardinals
(which is a weaker assumption than Vop\v{e}nka's principle) implies
that every extremally epireflective class $\Le$ is a small-orthogonality
class, under mild conditions on the category and the given class. These
conditions are fulfilled if the category is accessible \cite{AR} and
$\Le$ is defined by an absolute formula.

In order to explain the role played by absoluteness,
we note that, if one assumes that measurable cardinals exist, then
the reduction $P_{\Ae}$ mentioned above becomes the zero functor
in the category of groups, since if $\lambda$ is measurable then
$\Hom(\Z^{\lambda}/\Z^{<\lambda},\Z)\ne 0$ by \cite{DG}, so in fact
$P_{\Ae}\Z=0$ and therefore $P_{\Ae}$ kills all groups. 
Remarkably, this example shows
that one may ``define'' a functor $P_{\Ae}$, namely reduction with
respect to a certain class of groups, and it happens that the
conclusion of whether $P_{\Ae}$ is trivial or not depends on the
set-theoretical axioms adopted. Thus, such a functor is not
\textit{absolute\/} in the sense of model theory, that is, there is no
absolute formula in the usual language of set theory whose
satisfaction determines precisely $P_{\Ae}$ or its image. 
A formula (possibly containing parameters) is called absolute if, 
whenever it is satisfied in an inner model 
of set theory, it is also satisfied in the universe $V$ of all sets.
For instance, the class of modules over a ring $R$
is defined by an absolute formula with $R$ as a parameter.
On the other hand, statements involving cardinals, unbounded quantifiers or
choices may fail to be absolute.

We thank J.~Rosick\'y for his interest in this article and for showing us an example, described in Section~5,
of an epireflective class of graphs which is not a small-orthogonality class under the negation 
of Vop\v{e}nka's principle, even if supercompact cardinals are assumed to exist.  
This is another instance of a class that cannot be defined by any absolute formula.

Analogous situations occur in other areas of Mathematics. 
For example, if there exists a supercompact cardinal, then
all sets of real numbers that are definable by formulas 
whose quantifiers range only over real numbers and ordinals, 
and have only real numbers and ordinals 
as parameters, are Lebesgue measurable~\cite{SW}. In fact, in order
to prove the existence of nonmeasurable sets of real numbers, 
one needs to use the axiom of choice, a device that produces 
nondefinable objects~\cite{S}.

\section{Preliminaries from category theory}
\label{section2}
To make the paper readable for both category
theorists and set theorists, we will first recall a few basic
notions and facts from both fields. Classes that are not sets will
be called \textit{proper classes}.

A \textit{category\/} $\Ce$ consists of a (possibly proper) class
of \textit{objects\/} and pairwise disjoint sets $\Ce(X,Y)$, called
\textit{hom-sets}, for all objects $X$ and~$Y$, whose members are called
\textit{morphisms\/} from $X$ to~$Y$, together with 
associative composition functions
\[
\Ce(X,Y)\times\Ce(Y,Z)\longrightarrow\Ce(X,Z)
\]
for all $X$, $Y$, $Z$, and a distinguished
element $\identity_X\in\Ce(X,X)$ for all~$X$, which is a
unit for composition. A morphism is an \textit{isomorphism\/}
if it has a two-sided inverse. If $\Ce$ is a category, the notation
$X\in\Ce$ means that $X$ is an object of~$\Ce$.

A morphism $m\colon X\to Y$ is a \textit{monomorphism\/} if whenever two morphisms
$f$ and $g$ from an object $A$ to $X$ are given with $m\circ f = m\circ g$, 
the equality \hbox{$f=g$} follows. \textit{Epimorphisms\/} are defined dually.
A category is called \textit{balanced\/} if every morphism that is
both a monomorphism and an epimorphism is an isomorphism.
The category of rings and the category of graphs are important examples of nonbalanced categories. 
In this article, as in~\cite{AR}, a \textit{graph\/} will be a set~$X$ equipped
with a binary relation, where the elements of $X$ are called vertices and there is a
directed edge from $x$ to $y$ if and only if the pair $(x,y)$ is in the binary relation.
Each map of graphs is determined by the images of the vertices.
Hence, the monomorphisms of graphs are the injective maps,
and epimorphisms of graphs are maps that are surjective on vertices
(but not necessarily surjective on edges).

A monomorphism $m\colon X\to Y$ is \textit{strong\/} if, given any commutative square
\[
\xymatrix{
 P\ar[d]_{u}\ar[r]^{e} & Q\ar[d]^{v} \\
 X\ar[r]^{m} & Y
}
\label{strong}
\]
in which $e$ is an epimorphism, there is a unique morphism 
$f\colon Q\to X$ such that $f\circ e=u$ and $m\circ f=v$.
A monomorphism $m$ is \textit{extremal\/} if, whenever it factors as
$m=v\circ e$ where $e$ is an epimorphism, it follows that $e$ is an isomorphism. 
Split monomorphisms are strong, and strong monomorphisms are extremal.
If a morphism is both an extremal monomorphism and an epimorphism,
then it is necessarily an isomorphism, and,
if $\Ce$ is balanced, then all monomorphisms are extremal.
The dual definitions and similar comments apply to epimorphisms.

A \textit{subobject\/} of an object $X$ in a category $\Ce$ is an equivalence class of
monomorphisms $A\to X$, where $m\colon A\to X$ and $m'\colon A'\to X$
are declared equivalent if there are morphisms $u\colon A\to A'$ and
$v\colon A'\to A$ such that $m=m'\circ u$ and $m'=m\circ v$. For simplicity,
when we refer to a subobject $A$ of an
object~$X$, we view $A$ as an object equipped with a
monomorphism $A\to X$. A subobject is called strong (or extremal)
if the corresponding monomorphism is strong (or extremal).
The notion of a \textit{quotient\/} of an
object $X$ is defined, dually, as an equivalence class of
epimorphisms $X\to B$, under the corresponding equivalence relation.
A category is called \textit{well-powered\/} if~the subobjects 
of every object form a set, and it is called \textit{co-well-powered\/} 
if the quotients of every object form a set. 

A \textit{functor\/} $F$ from a category $\Ce$ to a category $\De$
associates to each object $X$ in $\Ce$ an object $FX$ in $\De$, and
to each morphism $f\colon X\to Y$ in $\Ce$ a morphism $Ff\colon
FX\to FY$ in~$\De$, preserving composition and identities. A functor
$F$ is \textit{full\/} if the function $\Ce(X,Y)\to\De(FX,FY)$ that sends
each morphism $f$ to $Ff$ is surjective for all $X$ and~$Y$, and it
is called \textit{faithful\/} if this function $\Ce(X,Y)\to\De(FX,FY)$
is injective for all $X$ and~$Y$. A subcategory $\Ae$ of a category
$\Ce$ is \textit{full\/} if the inclusion functor $\Ae\to\Ce$ is full.

A \textit{concrete\/} category is a category $\Ce$ together with a
faithful functor to the category of sets, $U\colon\Ce\to\Sets$.
See \cite{AHS} for an extensive treatment of this notion.
For an object $X$ of~$\Ce$, the set $UX$ is called the \textit{underlying set\/}
of~$X$, and similarly for morphisms.
In this article, when we assume that a category is concrete,
the functor $U$ will, as customary, 
be omitted from the notation. Hence we denote indistinctly an object $X$ of
$\Ce$ and its underlying set, and morphisms $X\to Y$ are also
seen as functions between the corresponding underlying sets.
In a concrete category, every morphism whose underlying
function is injective is a monomorphism, and every morphism whose
underlying function is surjective is an epimorphism.
Hence, for example, the homotopy category of topological spaces 
cannot be made concrete.

If $F$ and $G$ are functors from a category $\Ce$ to a category~$\De$,
a \textit{natural transformation\/} $\eta$ from $F$ to $G$ associates to
every object $X$ in $\Ce$ a morphism $\eta_X \colon FX \to GX$ in $\De$
such that, for every morphism $f \colon X \to Y$ in~$\Ce$, the
following diagram commutes:
\[
\xymatrix{
 FX\ar[d]_{\eta_X}\ar[r]^{Ff} & FY\ar[d]^{\eta_Y} \\
 GX\ar[r]^{Gf}                & GY.
}
\]

A \textit{reflection\/} (also called a \textit{localization\/}) 
on a category $\Ce$ is a pair $(L,\eta)$
where $L\colon\Ce\to\Ce$ is a functor and $\eta \colon \Id\to L$ is
a natural transformation, called \textit{unit}, 
such that $\eta_{LX}\colon LX\to LLX$ is an
isomorphism and $\eta_{LX}=L\eta_X$ for all $X$ in~$\Ce$. By abuse
of terminology, we often say that the functor $L$ itself is a reflection, 
or a localization, if the natural transformation $\eta$ is clear from the context.

If $L$ is a reflection, the objects $X$ such that
$\eta_X\colon X\to LX$ is an isomorphism
are called \textit{$L$-local objects}, and the morphisms $f$ such that $Lf$ is
an isomorphism are called \textit{$L$-equivalences}. By definition,
$\eta_X$ is an $L$-equivalence for all~$X$.
In fact, $\eta_X$ is terminal among $L$-equivalences with domain~$X$,
and it is initial among morphisms from $X$ to $L$-local objects.
This means that for each $L$-equivalence $h\colon X\to Y$ there is a unique
$h'\colon Y\to LX$ with $h'\circ h=\eta_X$, and
for each morphism $g\colon X\to Z$ where $Z$ is
$L$-local there is a unique $g'\colon LX\to Z$ such that $g'\circ\eta_X=g$.

A morphism $f\colon A\to B$ and an object $X$ are called \textit{orthogonal\/}
if for every morphism $g\colon A\to X$ there is a unique $g'\colon
B\to X$ such that $g'\circ f=g$. If $L$ is any reflection, then
an object is $L$-local if and only if it is orthogonal to all $L$-equivalences,
and a morphism is an $L$-equivalence if and only if it is orthogonal
to all $L$-local objects. See~\cite{CKort} for proofs or references
of these and other features of reflections.

A reflection $L$ is called an \textit{epireflection\/} if, for every $X$
in~$\Ce$, the unit morphism $\eta_X\colon X\to LX$ is an epimorphism.
We say that $L$ is a strong (or extremal) epireflection if $\eta_X$
is a strong (or extremal) epimorphism for all~$X$. 
A typical example of an epireflection is the abelianization functor
on the category of groups, associating to each group $G$ the quotient
by its commutator subgroup, equipped with the projection
$\eta_G\colon G\to G/[G,G]$.
The commutator subgroup is an example of a \textit{radical\/} on groups.
In the category of groups, there is a bijective correspondence between
epireflections and radicals, as shown in~\cite{CRS}.

Since a full subcategory is completely determined by the class of
its objects, the terms \textit{reflective class\/} and \textit{reflective
full subcategory\/} are used indistinctly to denote the class of $L$-local
objects for a reflection~$L$ or the full subcategory with these objects.
If $L$ is an epireflection, then the class of its local objects
is called \textit{epireflective}. It is called strongly epireflective
or extremally epireflective if $L$ is a strong or extremal epireflection.

The facts stated in the next proposition are not new.
Versions of them can be found in \cite{AHS} or~\cite{CRS}.

\begin{proposition}
\label{epiref}
Let $(L,\eta)$ be a reflection on a category~$\Ce$.
\begin{itemize}
\item[{\rm (a)}] If $L$ is an epireflection, then the class of $L$-local objects
is closed under strong subobjects, and it is closed under all subobjects if $\Ce$ is balanced.
\item[{\rm (b)}] Suppose that $\eta_X\colon X\to LX$ can be factored as an epimorphism followed
by a monomorphism for all~$X$. If the class of $L$-local objects is closed under subobjects, 
then $L$ is an epireflection. 
\end{itemize}
\end{proposition}

\begin{proof}
In order to prove~(a), let $s\colon A\to X$ be a monomorphism where
$X$ is $L$-local. By definition, $\eta_X$ is an isomorphism and hence the
composite $\eta_X\circ s$ is a monomorphism. Since $\eta_X\circ s=
Ls\circ\eta_A$ (because $\eta$ is a natural transformation), we
infer that $\eta_A$ is a monomorphism. Now, if $\Ce$ is balanced, then
$\eta_A$ is an isomorphism, so $A$ is $L$-local. If we assume instead
that $s$ is a strong monomorphism, then the existence of an inverse
of $\eta_A$ follows too.

To prove~(b), let $X$ be any object. Factor $\eta_X$ as
\[
X\stackrel{e}{\longrightarrow} Y\stackrel{m}{\longrightarrow} LX
\]
where $e$ is an epimorphism and $m$
is a monomorphism. Then $Y$ is a subobject of $LX$
and hence, by assumption, it is $L$-local.
Hence there is a unique morphism $f\colon LX\to Y$
such that $f\circ \eta_X=e$.
Then \[ m\circ f\circ\eta_X=m\circ e=\eta_X, \]
from which we infer that $m\circ f$ is the identity
(by the universal property of~$\eta_X$).
Hence $m$ is a split epimorphism and a monomorphism, from
which it follows that $m$ is an isomorphism.
\end{proof}

Note that, in part~(b), the conclusion that $L$ is an epireflection
also follows if ``monomorphism'' is replaced by strong (or extremal) monomorphism, 
and ``subobject'' is replaced by strong (or extremal) subobject. 
On the other hand, if ``epimorphism'' is replaced by strong (or extremal) 
epimorphism, then the argument used in the proof of part~(b) 
shows that $L$ is a strong (or extremal) epireflection.

A category is \textit{complete\/} if all set-indexed limits exist,
and it is \textit{cocomplete\/} if all set-indexed colimits exist.
See~\cite{AHS} or \cite{MacLane} for more information about limits and colimits,
and about products and coproducts in particular.

\newpage

\begin{proposition}
\label{epiexist}
If a category $\Ce$ is complete, well-powered, and co-well-powered, then
every class of objects $\Le$ closed under products and extremal subobjects
in $\Ce$ is epireflective, and if $\Le$ is closed under products and subobjects
then it is extremally epireflective.
\end{proposition}

\begin{proof}
It follows from~\cite[Proposition~12.5 and Corollary~14.21]{AHS} that, 
if $\Ce$ is complete and well-powered, then
every morphism in $\Ce$ can be factored as an extremal epimorphism followed by
a monomorphism, and also as an epimorphism followed by an extremal monomorphism. 
Thus, we may define a reflection by factoring, for each object~$X$, the
canonical morphism from $X$ into the product of its quotients that are 
in $\Le$ as an epimorphism $\eta_X$ followed by an extremal monomorphism,
or alternatively as an extremal epimorphism followed by a~monomorphism if $\Le$
is closed under subobjects.
\end{proof}

For each reflection $L$ on a category~$\Ce$, the class of $L$-local objects is closed
under all limits that exist in~$\Ce$, and the class of $L$-equivalences
is closed under all colimits that exist in the category of arrows of~$\Ce$
(whose objects are the morphisms of~$\Ce$ and whose morphisms are
commutative squares). In particular, every coproduct of $L$-equivalences
is an $L$-equivalence. 
If $\{ f_{i}\colon P_{i}\to Q_{i}\mid i\in I \}$ is a family of morphisms in
$\Ce$ and the coproducts $\coprod_{i\in I}P_{i}$ and
$\coprod_{i\in I}Q_{i}$ exist, with associated morphisms
$p_{i}\colon P_i\to \coprod_{i\in I}P_{i}$ and $q_{i}\colon Q_i\to
\coprod_{i\in I}Q_{i}$, then the coproduct $\coprod_{i\in I} f_i$
exists; namely, it is the unique morphism
\[
f \colon \coprod_{i\in I}P_{i}\longrightarrow \coprod_{i\in I}Q_{i}
\]
such that $f \circ p_{i}=q_{i}\circ f_i$ for all $i\in I$.

A \textit{small-orthogonality class\/} in a category $\Ce$ is the class of objects
orthogonal to some set of morphisms $\Fe=\{f_i\colon P_i\to Q_i\mid i\in I\}$.
An object orthogonal to all the morphisms in $\Fe$ will be called \textit{$\Fe$-local}.
If a reflection $L$ exists such that the class of $L$-local objects coincides with the 
class of $\Fe$-local objects for some set
of morphisms~$\Fe$, then $L$ will be called an \textit{$\Fe$-localization\/}
(or an \textit{$f$-localization\/} if $\Fe$ consists of one morphism $f$ only).

Note that, if a coproduct $f=\coprod_{i\in I} f_i$ exists and all
hom-sets $\Ce(X,Y)$ of $\Ce$ are nonempty, then an object is orthogonal
to $f$ if and only if it is orthogonal to $f_i$ for all $i\in I$. 
More precisely, if $X$ is orthogonal to all $f_i$ then it is orthogonal to their coproduct,
and the converse holds if $\Ce(P_i,X)\ne\emptyset$ for all $i\in I$, where $P_i$
is the domain of~$f_i$.
Hence, if $\Ce$ has coproducts and all its hom-sets are nonempty, then every 
small-orthogonality class is the class of objects orthogonal to a single morphism.

A sufficient condition for a category ensuring that all hom-sets are non\-empty
is the existence of a \textit{zero object}, that is, an object $0$ which
is both initial and final. This is the case, for example, with the
trivial group in the category of groups and with the one-point
space in the category of topological spaces with a base point. If
$\Ce$ has a zero object, then each set $\Ce(X,Y)$ contains at least the
\textit{zero morphism\/} $X\to 0\to Y$. 

\newpage

\begin{proposition}
\label{epis}
Let $(L,\eta)$ be an $\Fe$-localization on a category $\Ce$, where $\Fe$ is a 
nonempty set of morphisms.
\begin{itemize}
\item[(a)] Suppose that every morphism of $\Ce$ can be factored
as an epimorphism followed by a strong monomorphism.
If every $f\in\Fe$ is an epimorphism, then $L$ is an epireflection.
\item[(b)] If $L$ is an epireflection, then there is a set $\Ee$
of epimorphisms such that $L$ is also an $\Ee$-localization.
\end{itemize}
\end{proposition}

\begin{proof}
By Proposition~\ref{epiref} (and the remark after it), 
in order to prove~(a) it suffices to check that the
class of $L$-local objects is closed under strong subobjects. Thus, let $X$ be $L$-local 
and let $s\colon A\to X$ be a strong monomorphism. We need to show that $A$
is orthogonal to every morphism $f\colon P\to Q$ in~$\Fe$.
For this, let $g\colon P\to A$ be any morphism. Since $X$ is orthogonal to~$f$, there is a
unique morphism $g'\colon Q\to X$ such that $g'\circ f=s\circ g$.
Since $f$ is an epimorphism and $s$ is strong, there is a morphism $g''\colon
Q\to A$ such that $g''\circ f= g$ and $s\circ g''=g'$.
Moreover, if $g'''\colon Q\to A$ also satisfies $g'''\circ f=g$,
then $g'''=g''$ since $f$ is an epimorphism. Hence, $A$ is orthogonal to~$f$.

Our argument for part~(b) is based on a similar result in~\cite{RS}. Write
$\Fe=\{f_i\colon P_i\to Q_i\mid i\in I\}$, and let
\[ \Ee=\{\eta_{P_i}\colon P_i\to LP_i\mid i\in I\} \cup \{\eta_{Q_i}\colon Q_i\to LQ_i\mid i\in I\}. \]
Then every morphism in $\Ee$ is an epimorphism, and the class of $\Ee$-local objects coincides precisely 
with the class of $\Fe$-local objects.
\end{proof}

\begin{example}
{\rm
In the category of graphs, let $L$
be the functor assigning to every graph $X$ the complete graph (i.e., containing
all possible edges between its vertices) with the same
set of vertices as~$X$, and let $\eta_X\colon X\to LX$ be the inclusion.
Then $L$ is an epireflection. The class of $L$-local objects is the class of
complete graphs, which is closed under strong subobjects, but not under
arbitrary subobjects. In fact $L$ is an $f$-localization, where $f$ is the
inclusion of the two-point graph $\{0,1\}$ into $0\to 1$, which is an epimorphism.
}
\end{example}

We finally recall the definition of locally presentable and accessible categories.
For a regular cardinal $\lambda$, a partially ordered set is called 
\textit{$\lambda$-directed\/} if every subset of cardinality smaller than
$\lambda$ has an upper bound. An object $X$ of a category $\Ce$ is called
\textit{$\lambda$-presentable}, where $\lambda$ is a regular cardinal, if
the functor $\Ce(X,-)$ preserves $\lambda$-directed colimits,
that is, colimits of diagrams indexed by $\lambda$-directed partially ordered sets. 
A category $\Ce$ is \textit{locally presentable\/} if it is
cocomplete and there is a regular cardinal $\lambda$ and a set
$\Xe$ of $\lambda$-presentable objects such that every object of
$\Ce$ is a $\lambda$-directed colimit of objects from~$\Xe$.
Locally presentable categories are complete, well-powered and co-well-powered.
The categories of groups, rings, modules over
a ring, and many others are locally presentable; see~\cite[1.B]{AR}
for further details and more examples.

If the assumption of cocompleteness is weakened by imposing instead that
$\lambda$-directed colimits exist in~$\Ce$, then $\Ce$ is
called \textit{$\lambda$-accessible}. A category $\Ce$ is called
\textit{accessible\/} if it is $\lambda$-accessible for some regular
cardinal~$\lambda$. As shown in~\cite[Theorem~5.35]{AR},
the accessible categories are precisely the categories
equivalent to categories of models of basic theories.
The definition of the latter terms is recalled at the end of
the next section.

\section{Preliminaries from set theory}
\label{section3}

The \textit{universe\/} $V$ of all sets is a proper class defined recursively on
the class $\OR$ of ordinals as follows: $V_0=\emptyset$,
$V_{\alpha +1}={\mathcal P}(V_{\alpha})$ for all~$\alpha$, 
where ${\mathcal{P}}$ is the power-set operation, and $V_{\lambda}=\bigcup_{\alpha
<\lambda}V_{\alpha}$ if $\lambda$ is a limit ordinal. Finally, $V=\cup_{\alpha\in\OR}V_{\alpha}$. 
Transfinite induction shows that, if $\alpha$ is any ordinal, then $\alpha\subseteq V_{\alpha}$.
The axiom of regularity, stating that every nonempty set has a
minimal element with respect to the membership relation, implies
that every set is an element of some $V_{\alpha}$; see~\cite[Lemma 9.3]{Jech1}.
The \textit{rank\/} of a set $X$, denoted $\rank(X)$, is the least
ordinal $\alpha$ such that $X\in V_{\alpha +1}$. Thus,
$\rank(\alpha)=\alpha$ for all ordinals $\alpha$.
More generally, if $X$ is any set, then $\rank(X)$ is the 
supremum of the set $\{ \rank(x)+1\mid x\in X\}$. 

A set or a proper class $X$ is called \textit{transitive\/} if every
element of an element of $X$ is also an element of~$X$. The universe $V$
is transitive, and so is $V_{\alpha}$ for every ordinal $\alpha$.
The \textit{transitive closure\/} of a set $X$, written $\TC(X)$,
is the smallest transitive set containing~$X$, that is, the intersection
of all transitive sets that contain~$X$. The elements of $\TC(X)$
are the elements of~$X$, the elements of the elements of~$X$, etc.

The \emph{language of set theory\/} is the first-order language
whose only nonlogical symbols are equality $=$ and the binary relation 
symbol~$\in$. The language consists of \emph{formulas\/} built up from the
\emph{atomic formulas\/} $x = y$ and $x \in y$, where $x$ and $y$ are
members of a set of variables, using the logical connectives $\neg$, $\wedge$, $\vee$,
$\to$, $\leftrightarrow$, and the quantifiers $\forall v$ and
$\exists v$, where $v$ is a variable. We use Greek letters to denote formulas.
The variables that appear in a formula $\varphi$ outside the scope of a quantifier are
called \emph{free}. The notation $\varphi (x_1,\ldots,x_n)$ means that 
$x_1,\ldots,x_n$ are the free variables in~$\varphi$.

All axioms of ZFC can be formalized in the language of set theory. 
A \emph{model\/} of ZFC is a set or a proper class $M$ in which the
formalized axioms of ZFC are true when the binary relation symbol
$\in$ is interpreted as the membership relation. A model $M$ is
called \emph{inner\/} if it is transitive and contains all
the ordinals. Thus, inner models are not sets, but proper classes.
Given a model $M$ and a formula $\varphi (x_1,\ldots,x_n)$, 
and given an $n$-tuple $a_1,\ldots,a_n$ of elements of~$M$, we say that
\emph{$\varphi (a_1,\ldots,a_n)$ is satisfied in $M$\/} if the formula
is true in $M$ when $x_i$ is replaced by $a_i$ for all $1\leq i\leq n$.

A set or a proper class $C$ is \textit{definable\/} in a
model $M$ if there is a formula $\varphi(x,x_1,\ldots,x_n)$ of the
language of set theory and elements $a_1,\ldots,a_n$ in $M$ such
that $C$ is the class of elements $c\in M$ such that
$\varphi(c,a_1,\ldots,a_n)$ is satisfied in~$M$.
We then say that $C$ is \textit{defined by $\varphi$ with parameters\/} $a_1,\ldots ,a_n$.
Notice that every set $a\in M$ is definable in $M$ with $a$ as a
parameter, namely by the formula $x\in a$.

A formula $\varphi (x,x_1,\ldots,x_n)$ is \textit{absolute between two
models\/} $N\subseteq M$ with respect to a collection of parameters $a_1,\ldots,a_n$ in $N$
if, for each $c\in N$, $\varphi(c,a_1,\ldots,a_n)$ is satisfied in $N$ if and only if it is satisfied
in~$M$. A formula is called \textit{absolute\/} with respect to $a_1,\ldots,a_n$ if it is absolute
between any inner model $M$ that contains $a_1,\ldots,a_n$
and the universe~$V$. We call a set or a proper class
\emph{absolute\/} if it is defined in $V$ by an absolute formula.

A submodel $N$ of a model $M$ is \textit{elementary\/}
if all formulas are absolute between $N$ and $M$ with respect
to every set of parameters in~$N$.
An embedding of $V$ into a model $M$ is an \textit{elementary embedding\/}
if its image is an elementary submodel of~$M$.
If $j\colon V\to M$ is a nontrivial elementary embedding with $M$ transitive,
then $M$ is inner, and induction on rank shows that
there is a least ordinal $\kappa$ moved by~$j$, that is, $j(\alpha)=\alpha$ for all $\alpha <\kappa$, 
and $j(\kappa)>\kappa$. Such a $\kappa$ is called the \textit{critical point\/} of~$j$,
and it is necessarily a measurable cardinal; see~\cite[Lemma~28.5]{Jech1}.

For a set $X$ and a cardinal~$\kappa$, let ${\mathcal P}_{\kappa}(X)$ be the set of subsets of $X$ 
of cardinality less than~$\kappa$. A cardinal $\kappa$ is called
\textit{$\lambda$-supercompact}, where $\lambda$ is an ordinal, if
the set ${\mathcal P}_{\kappa}(\lambda)$ admits a normal measure~\cite{Jech1}.
A cardinal $\kappa$ is \textit{supercompact\/} if it is
$\lambda$-supercompact for every ordinal~$\lambda$. 
Instead of recalling the definition of a normal measure,
we recall from \cite[Lemma~33.9]{Jech1} that a cardinal $\kappa$ is
$\lambda$-supercompact if and only if there is an elementary
embedding $j\colon V\to M$ such that $j(\alpha)=\alpha$ for all
$\alpha<\kappa$ and $j(\kappa)>\lambda$, where $M$ is an inner model
such that $\{ f\mid f\colon \lambda \to M\} \subseteq M$,
i.e., every $\lambda$-sequence of elements of $M$ is an element of~$M$. 
For more information on supercompact cardinals, see \cite{Jech2} or~\cite{Kanamori}. 

If $j\colon V\to M$ is an elementary embedding, then 
for every set $X$ the \emph{restriction\/} $j\restriction X\colon X\to j(X)$ is the
function that sends each element $x\in X$ to~$j(x)$.
The statement that $j\restriction X\colon X\to j(X)$ is in $M$
means that the set $\{(x,j(x))\mid x\in X\}$ is an element of~$M$.

\begin{proposition}
A cardinal $\kappa$ is supercompact if and only if for every set $X$
there is an elementary embedding $j$ of the universe $V$ into an inner model $M$
with critical point~$\kappa$, such that $X\in M$, $j(\kappa)>\rank(X)$, 
and $j\restriction X\colon X\to j(X)$ is in~$M$.
\end{proposition}

\begin{proof}
Given any set~$X$, let $\lambda$ be the cardinality of the transitive
closure of the set~$\{X\}$, and consider the binary relation $R$ on $\lambda$
that corresponds to the membership relation on this transitive closure.
By~\cite[3.12]{Jech2}, the binary relation $R$ embeds into~$\lambda$.
Therefore, the set $X$ is encoded by a \hbox{$\lambda$-sequence} of ordinals.
Now choose an elementary embedding $j\colon V\to M$
with $M$ transitive and critical point~$\kappa$, such that $j(\kappa)>\lambda$
and $M$ contains all the \hbox{$\lambda$-sequences} of its elements. From
the latter it follows that $X\in M$.
Finally, we use the fact that the restriction $j\restriction\lambda$ is in $M$
if and only if $\{ f\mid f\colon \lambda \to M\} \subseteq M$;
see~\cite[Proposition~22.4]{Kanamori}.
\end{proof}

We finally recall the following definitions from~\cite[Chapter~5]{AR}.
For a set~$S$ and a regular cardinal~$\lambda$, a \textit{$\lambda$-ary $S$-sorted signature\/} 
$\Sigma$ consists of a set of \textit{operation symbols}, each of which has a certain
\textit{arity\/} $\prod_{i\in I} s_i\to s$, where $s$ and all $s_i$ are in~$S$
and $|I|<\lambda$, and another set of \textit{relation symbols}, each of which has also a certain arity
of the form $\prod_{j\in J} s_j$, where all $s_j$ are in $S$ and $|J|<\lambda$.
Given a signature~$\Sigma$, a \textit{$\Sigma$-structure\/} is a collection
$X=\{X_s\mid s\in S\}$ of nonempty sets together with a function
\[ \sigma_X\colon \prod_{i\in I}X_{s_i}\longrightarrow X_s \]
for each operation symbol $\sigma\colon\prod_{i\in I}s_i\to s$, and
a subset $\rho_X\subseteq \prod_{j\in J}X_{s_j}$ for each relation symbol $\rho$
of arity $\prod_{j\in J}s_j$. A \textit{homomorphism\/} of
$\Sigma$-structures is a collection $f=\{f_s\mid s\in S\}$ of functions
preserving operations and relations. The category of $\Sigma$-structures and
their homomorphisms is denoted by~${\bf Str}\,\Sigma$.

Given a $\lambda$-ary $S$-sorted signature $\Sigma$
and a collection $W=\{W_s\mid s\in S\}$ of sets of cardinality~$\lambda$, where the elements
of $W_s$ are called \textit{variables of sort~$s$}, one defines \textit{terms\/} by declaring
that each variable is a term and, for each operation symbol
$\sigma\colon\prod_{i\in I}s_i\to s$ and each collection of terms $\tau_i$
of sort~$s_i$, the expression $\sigma(\tau_i)_{i\in I}$
is a term of sort~$s$. \textit{Formulas\/} are built up by means of 
logical connectives and quantifiers from the \textit{atomic formulas\/} $\tau_1=\tau_2$ and
$\rho(\tau_j)_{j\in J}$, where $\rho$ is a relation symbol and each $\tau_j$ is a term.
Variables which appear unquantified in a formula are said to appear free. 
A formula without free variables is called
a \textit{sentence}. A set of sentences is called a \textit{theory\/} (with signature $\Sigma$).
A \textit{model\/} of a theory $T$ with signature $\Sigma$ is a $\Sigma$-structure satisfying each sentence of~$T$.
For each theory~$T$, we denote by ${\bf Mod}\,T$ the full subcategory of
${\bf Str}\,\Sigma$ consisting of all models of~$T$.

A formula is called \textit{basic\/} if it has the form $\forall x(\varphi(x)\to
\psi(x))$, where $\varphi$ and $\psi$ are disjunctions of formulas of type
$\exists y\;\zeta(x,y)$ in which $\zeta$ is a conjunction of atomic formulas.
A \textit{basic theory\/} is a theory of basic sentences. By ~\cite[Theorem~5.35]{AR}, 
a category is accessible if and only if it is equivalent to ${\bf Mod}\,T$ for some basic theory~$T$.

\section{Main results}
\label{section4}

If $\Ae$ is a class of objects in a category $\Ce$,
a set $\He$ of objects of $\Ce$ will be called \textit{transverse\/} to
$\Ae$ if every object of $\Ae$ has a subobject in~$\He\cap\Ae$.

\begin{theorem}
\label{firsttheorem} 
Suppose that $(L,\eta)$ is an epireflection on a category $\Ce$.
\begin{itemize}
\item[(a)] If $\Ce$ is balanced and there exists a
set $\He$ of objects in $\Ce$ which is 
transverse to the class of objects that are not $L$-local,
then there is a set of morphisms $\Fe$ such that $L$ is an $\Fe$-localization. 
\item[(b)] If $\Ce$ is co-well-powered and every morphism can be factored
as an epimorphism followed by a monomorphism,
then the converse holds, that is, if $L$ is an
$\Fe$-localization for some set of morphisms~$\Fe$, then there is a 
set $\He$ transverse to the class of objects that are not $L$-local.
\end{itemize}
\end{theorem}

\begin{proof}
To prove~(a), let $\Fe=\{\eta_A\colon A\to LA\mid A\in\He\}$.
Fix any object $X$ of~$\Ce$. If $X$ is $L$-local, then
$X$ is orthogonal to all morphisms in $\Fe$, since these are $L$-equivalences.
In other words, $X$ is $\Fe$-local.
Now suppose that $X$ is not \hbox{$L$-local}. We aim to show that $X$
is not $\Fe$-local, hence completing the proof. By assumption, in the
set $\He$ there is a subobject $A$ of $X$ that is not \hbox{$L$-local}.
Let $s\colon A\to X$ be a monomorphism. 
Towards a contradiction, suppose that $X$ is $\Fe$-local. Then $X$ is
orthogonal to~$\eta_{A}$. Hence there is a morphism $t\colon LA\to
X$ such that $s=t\circ \eta_{A}$. This implies that $\eta_{A}$ is a
monomorphism and hence an isomorphism, since $\Ce$ is balanced. This
contradicts the fact that $A$ is not isomorphic to~$LA$. Hence,
$X$ is not $\Fe$-local, as needed.

For the converse, suppose that $L$ is an $\Fe$-localization for some
nonempty set of morphisms $\Fe=\{f_i\colon P_i\to Q_i\mid i\in I\}$. 
Since $L$ is an epireflection, we may assume, by part~(b) of Proposition~\ref{epis},
that each $f_i$ is an epimorphism. Since we suppose that $\Ce$ is co-well-powered, 
we may consider the set $\He$ of all quotients of $P_i$ for all $i\in I$ (that is, we choose a
representative object of each quotient). Let $X$ be an object which is not
$L$-local. Note that, if a morphism $P_i\to X$ can be factored
through $Q_i$, then it can be factored in a unique way, since $f_i$ is an
epimorphism. Hence, if $X$ is not $L$-local, then there is a morphism 
$g\colon P_i\to X$ for some $i\in I$ for which there is
no morphism $h\colon Q_i\to X$ with $h\circ f_i=g$. Factor $g$ as
$g''\circ g'$, where $g'\colon P_i\to X'$ is an epimorphism and
$g''\colon X'\to X$ is a monomorphism, in such a way that $X'$ is in~$\He$.
Note finally that $X'$ is not $L$-local, for if it were then there would exist a morphism
$h'\colon Q_i\to X'$ such that $g''\circ h'\circ f_i=g$, which, as we know, cannot happen.
\end{proof}

\begin{remark}
\label{notbalanced}
{\rm
For the validity of part~(a) of Theorem~\ref{firsttheorem}, 
the assumption that $\Ce$ is balanced can
be weakened by assuming only that the epimorphisms $\eta_A$ are
extremal for $A\in\He$, so that they are isomorphisms whenever they are monomorphisms. 
This ensures the validity of the theorem in important categories that are not
balanced, such as the category of graphs (see Section~5 below), provided that
$L$ is an extremal epireflection.
By Proposition~\ref{epiref}, the condition that $L$ is an extremal
epireflection is satisfied if the class of $L$-local objects
is closed under subobjects, and morphisms in $\Ce$ can be factored as an extremal
epimorphism followed by a monomorphism. By~\cite[Proposition~1.61]{AR},
the latter holds in locally presentable categories. More generally,
it holds in complete well-powered categories, by~\cite[Corollary~14.21]{AHS}.

Note also that, if we add the assumption that 
$\Ce$ has coproducts and $\Ce(X,Y)$ is nonempty for all objects $X$ and $Y$, then the set of
morphisms $\Fe$ given by part (a) of the theorem can be replaced by a single morphism~$f$,
namely the coproduct of all morphisms in~$\Fe$.
}
\end{remark}

In the rest of this section, all categories will be assumed to be concrete,
and the corresponding faithful functor into the category of sets will be omitted from the notation.
A concrete category $\Ce$ will be called \textit{absolute\/} if there is an absolute formula $\varphi
(x,y,z,x_1,\ldots,x_n)$ with respect to a set of parameters $a_1,\ldots,a_n$ such that, for any two sets $A$,
$B$ and any function $f\colon A\to B$, $\varphi (A,B,f,a_1,\ldots,a_n)$ is satisfied
in the universe $V$ if and only if $A$ and $B$ are objects of $\Ce$ and $f$ is in~$\Ce(A,B)$. 
For example, the categories of groups, rings, or modules over a ring $R$ are absolute.
(In the latter case, the ring $R$ is a parameter; in the other two examples, there are no parameters.)
More generally, every category ${\bf Mod}\, T$ of models over a theory $T$ is absolute.
Therefore, by~\cite[Theorem~5.35]{AR}, all accessible categories are absolute.

A reflection $L$ will be called \textit{absolute\/} if the class of $L$-local objects is absolute.
For example, abelianization of groups is absolute, and, more generally,
every projection onto a variety of groups is absolute; see~\cite{CRS}.

\begin{definition}
{\rm We say that a concrete category $\Ce$ \textit{supports elementary
embeddings\/} if, for every elementary embedding $j\colon V\to M$ 
and all objects $X$ of~$\Ce$, the restriction
$j\restriction X\colon X\to j(X)$ underlies a morphism of~$\Ce$.}
\end{definition}

Note that $j\restriction X\colon X\to j(X)$ is always injective,
since $j(x)=j(y)$ implies that $x=y$.
Hence, if $\Ce$ is concrete and supports elementary embeddings, then
$j\restriction X$ is a monomorphism in $\Ce$ for all~$X$.

\begin{proposition}
If $\Ce$ is an absolute full subcategory of ${\bf Str}\,\Sigma$ for some
signature~$\Sigma$, then $\Ce$ supports elementary embeddings.
\end{proposition}

\begin{proof}
We first prove that ${\bf Str}\,\Sigma$ itself supports elementary
embeddings. If $X$ is a $\Sigma$-structure, then the set $j(X)$ admits
operations and relations defined as $\sigma_{j(X)}=j(\sigma_X)$
for every operation symbol $\sigma$ of~$\Sigma$, and $\rho_{j(X)}=j(\rho_X)$
for every relation symbol~$\rho$. Thus, $j(X)$ becomes a $\Sigma$-structure
in such a way that $j\restriction X\colon X\to j(X)$ is a homomorphism
of $\Sigma$-structures.

Now let $\Ce$ be an absolute full subcategory of~${\bf Str}\,\Sigma$.
If $X$ is an object in~$\Ce$ then~$j(X)$, viewed as a $\Sigma$-structure as in the previous
paragraph, is also an object of $\Ce$ since $\Ce$ is assumed to be absolute, 
and the function $j\restriction X$ is automatically a homomorphism of $\Sigma$-structures.
Since $\Ce$ is assumed to be full, $j\restriction X$ is a morphism in~$\Ce$.
\end{proof}

Therefore, by~\cite[Theorem~5.35]{AR}, accessible categories support elementary embeddings. 
Accessible categories are indeed concrete, since they
can be embedded into the category of graphs~\cite[Theorem~2.65]{AR}.

It is however not true that every absolute concrete category supports elementary embeddings.
For example, let $\Ce$ be the category whose class of objects is the class $V$ of all sets
and whose morphisms are defined by $\Ce(X,Y)=\emptyset$ if $X\ne Y$ and $\Ce(X,X)=\{{\rm id}_X\}$
for all~$X$. Then $\Ce$ does not support elementary embeddings.

\begin{theorem}
\label{theorem3.3} Suppose that $\kappa$ is a supercompact cardinal and $\Ae$
is an absolute class of objects in an absolute category $\Ce$
which supports elementary embeddings. Suppose also that the parameters in the 
definitions of $\Ce$ and $\Ae$ have rank less than~$\kappa$. If $X\in\Ae$,
then there is a subobject of $X$ in $V_{\kappa}\cap\Ae$.
\end{theorem}

\begin{proof}
Let $\varphi$ be an absolute formula defining $\Ae$ in $V$ with parameters
$a_1,\ldots,a_n$, and let $b_1,\ldots,b_m$ be the parameters in the definition
of the category~$\Ce$. Fix an object $X\in\Ae$
and let $j\colon V\to M$, with $M$ transitive, be
an elementary embedding with critical point $\kappa$ such that $X$
and the restriction $j\restriction X$ are~in~$M$, and
$j(\kappa)>\mbox{rank}(X)$. Notice that $a_1,\ldots,a_n$ and $b_1,\ldots,b_m$ are also in~$M$,
since in fact $j(a_r)=a_r$ for all~$r$ and $j(b_s)=b_s$ for all~$s$. 
Let us write $\vec{a}$ for $a_1,\ldots ,a_n$ and $\vec{b}$ for $b_1,\ldots ,b_m$.

Since $\Ce$ is absolute, $j(X)$ is an object of~$\Ce$. 
Moreover, since $\Ce$ supports elementary embeddings,
the restriction $j\restriction X\colon X\to j(X)$ underlies a
monomorphism in~$\Ce$. Hence, $j(X)$ has a subobject in~$M$,
namely~$X$, which satisfies $\varphi$ and has rank less
than~$j(\kappa)$. Now ``$y$ is a subobject of $x$'' means
``$x$~and $y$ are objects of $\Ce$ and there is a morphism $y\to x$ which is a monomorphism''.
Hence, the following formula in the parameters $X$, $\vec{a}$, $\vec{b}$, $\kappa$ 
is true in~$M$:
\[ \exists y\, ((y \mbox{ is a subobject of }j(X)) \wedge \varphi
(y, \vec{a}) \wedge ({\rm rank}(y)<j(\kappa))). \]
Hence, since $j$ is an elementary embedding, the following holds in~$V$:
\[ \exists y\, ((y \mbox{ is a subobject of }X) \wedge \varphi
(y,\vec{a}) \wedge ({\rm rank}(y)<\kappa)). \] That is, $X$ has a
subobject in $V_{\kappa}\cap\Ae$, which proves the theorem.
\end{proof}

\begin{corollary}
\label{corollary3.4} Suppose that $(L,\eta)$ is an absolute extremal
epireflection on an absolute category $\Ce$ which supports
elementary embeddings. If there is a supercompact
cardinal $\kappa$ greater than the ranks of the parameters in the definition of $\Ce$ and in the definition of 
the class of $L$-local objects, then $L$ is an \hbox{$\Fe$-localization} for some set $\Fe$ of morphisms.
\end{corollary}

\begin{proof}
Let the class of objects of $\Ce$ that are not $L$-local play the role of the
class $\Ae$ in Theorem~\ref{theorem3.3}. Then the conclusion of the theorem
is precisely that the set $V_{\kappa}$ is transverse to the class of objects
of $\Ce$ that are not $L$-local. Hence, part~(a) of Theorem~\ref{firsttheorem}
and Remark~\ref{notbalanced} yield the desired result.
\end{proof}

Recall that, if $\Ce$ is balanced, then every epireflection is extremal.
And if we assume that $\Ce$ has coproducts and $\Ce(X,Y)$ is nonempty for all $X$ and~$Y$,
then we may infer, in addition to the conclusion of Corollary~\ref{corollary3.4}, 
that $L$ is an $f$-localization for a single morphism~$f$,
which can be chosen to be an epimorphism by Proposition~\ref{epis}.

As an application, we give the following result.
For any given class of groups~$\Ae$,
the \textit{reduction\/} $P_{\Ae}$ is an epireflection on the
category of groups whose local objects are groups
$G$ that are \textit{$\Ae$-reduced}, i.e., for which 
every homomorphism $A\to G$ is trivial if $A\in\Ae$.
Such an epireflection exists by Proposition~\ref{epiexist}, since the class of
$\Ae$-reduced groups is closed under products and subgroups.

\begin{corollary}
\label{application}
Let $\Ae$ be any absolute class of groups (possibly proper).
If there is a supercompact cardinal greater than the ranks of the parameters 
in the definition of $\Ae$, then there is a group $U$
such that the class of $U$-reduced groups coincides with
the class of $\Ae$-reduced groups.
\end{corollary}

\begin{proof}
The category of groups is balanced and locally presentable.
Hence, Corollary~\ref{corollary3.4} implies that
the reduction functor $P_{\Ae}$ is an $f$-localization for
some group homomorphism~$f$. As in~\cite[Theorem~6.3]{CRS},
let $U$ be a universal $f$-acyclic group, i.e., a group $U$
such that $P_U$ and $P_{\Ae}$ annihilate the same groups.
Then, by \cite[Theorem~2.3]{CRS}, $P_U$ and $P_{\Ae}$
also have the same class of local objects; that is,
the class of $U$-reduced groups coincides indeed
with the class of $\Ae$-reduced groups.
\end{proof}

As pointed out in the Introduction, for the
(non-absolute) class $\Ae$ of groups of the form $\Z^{\kappa}/\Z^{<\kappa}$ for all
cardinals $\kappa$, the existence of a group $U$ such that the
class of $U$-reduced groups coincides with the class of $\Ae$-reduced groups
is equivalent to the existence of a measurable cardinal; see~\cite{CRS} or~\cite{DG}.

\section{On absoluteness}

We will display an example, indicated to us by Rosick\'y, 
of an extremal epireflection $L$ on the category $\Gra$ of graphs
which is not an $\Fe$-localization for any set of maps~$\Fe$.
This example is based on \cite[Example~6.12]{AR} and requires to
assume the negation of Vop\v{e}nka's principle while admitting
the existence of supercompact cardinals.

Since we are assuming that Vop\v{e}nka's principle does not hold, we may
choose a proper class of graphs $\Ae$ which is \textit{rigid}, that is, such that 
\[ \Gra(A,B)=\emptyset\]
for all $A\ne B$ in~$\Ae$, and $\Gra(A,A)$ has the identity as its only element 
for every $A\in\Ae$.
Consider the class $\Le$ of graphs that are \textit{$\Ae$-reduced}, i.e.,
\[ \Le=\{ X\in \Gra\mid\Gra(A,X)=\emptyset \mbox{ for all }A\in\Ae \}, \]
and note that $\Ae\cap\Le=\emptyset$, while every proper subgraph of a graph
in $\Ae$ is in~$\Le$. By Proposition~\ref{epiexist}, 
there is an epireflection $L$ whose class of local objects is precisely~$\Le$,
since $\Le$ is closed under products and subobjects in the category of graphs.
Moreover, the unit map $\eta_X\colon X\to LX$ is an extremal epimorphism 
(indeed, surjective on vertices and edges) for all~$X$.

Now suppose that there is a set $\Fe=\{f_i\colon P_i\to Q_i\mid i\in I\}$ of maps of graphs
such that the reflection $L$ is an $\Fe$-localization. Then, if we choose any
regular cardinal $\lambda$ that is bigger than the cardinalities of
$P_i$ and $Q_i$ for all $i\in I$, it follows that $\Le$ is closed under $\lambda$-directed colimits.
As in \cite[Example~6.12]{AR}, a contradiction is obtained by
choosing a graph $A\in\Ae$ whose cardinality is bigger than~$\lambda$,
and observing that $A$ is a $\lambda$-directed colimit of 
the diagram of all its proper subgraphs, each of which is in~$\Le$,
while $A$ itself is not in~$\Le$. This contradicts the
previous statement that $\Le$ is closed under $\lambda$-directed colimits.

Thus, we infer that the class $\Le$ cannot be absolute, since otherwise this example
would contradict Corollary~\ref{corollary3.4}.
The fact that $\Le$ is not absolute can be seen directly as follows.
Suppose that $\Le$ is absolute, so $\Ae$ is also absolute. 
Let $\varphi$ be a formula defining $\Ae$ (possibly with parameters) and $\psi$ the corresponding
formula defining $\Le$, namely 
\[ (x\in\Gra)\wedge \forall y(((y\in\Gra)\wedge\varphi(y))\to\Gra(y,x)=\emptyset) .\]
Let $\kappa$ be a supercompact cardinal and
choose a graph $A\in\Ae$ with $|A|>\kappa$. Let $\lambda$ be a regular cardinal
such that $\lambda >|A|$. Since $\kappa$ is supercompact, there is an elementary embedding
$j\colon V\to M$ with critical point $\kappa$ such that $j(\kappa)>\lambda$, $A\in M$, and
$j\restriction A\colon A\to j(A)$ is also in~$M$. 
Note that $j(A)$ is a graph, since elementary embeddings preserve binary relations. 
From the fact that $j\restriction A$ is in $M$ it follows that $A$ is a subgraph of $j(A)$ in~$M$,
and moreover it is proper subgraph, since
\[ |A|<\lambda<j(\kappa)<|j(A)|, \]
where the last inequality follows from the fact that $\kappa<|A|$.
Then $A$ satisfies the formula $\psi$ in~$M$, since it is a proper subgraph of a graph 
satisfying~$\varphi$, namely~$j(A)$.
Since $j$ is elementary, $A$ also satisfies $\psi$ in~$V$, that is, $A\in\Le$. 
Hence $A\in\Ae\cap\Le$, which contradicts the fact that $\Ae\cap\Le=\emptyset$.

This example shows in fact that, if there are supercompact cardinals, then Vop\v{e}nka's principle
holds for absolute classes of graphs defined with small parameters;
that is, for a supercompact cardinal~$\kappa$, there is no rigid absolute proper class of graphs 
defined with parameters of cardinality smaller than~$\kappa$.

\vskip 0.5cm

{\small

\noindent
Joan Bagaria,
ICREA (Instituci\'o Catalana de Recerca i Estudis Avan\c{c}ats) and
Departament de L\`ogica, Hist\`oria i Filosofia de la Ci\`encia, Universitat de Barcelona. \\
Montalegre 6,
08001 Barcelona, Spain. \\
{\tt joan.bagaria@icrea.es}

\bigskip

\noindent
Carles Casacuberta,
Departament d'\`Algebra i Geometria, Universitat de Barcelona. \\
Gran Via de les Corts Catalanes 585, 08007 Barcelona, Spain. \\
{\tt carles.casacuberta@ub.edu}

\bigskip

\noindent
Adrian R. D. Mathias,
ERMIT, Universit\'e de la R\'eunion. \\
Avenue Ren\'e Cassin 15,
F-97715 Saint Denis de la R\'eunion, France. \\
{\tt ardm@dpmms.cam.ac.uk}

}

\end{document}